\documentclass{elsart} 

\usepackage{color}
\usepackage[latin1]{inputenc}
\usepackage[dvips]{graphicx}
\usepackage{amssymb}

\newcommand{\melo}{Glaucio Gomes de Magalhães Melo}
\newcommand{\oliveiralima}{Emerson Alexandre de Oliveira Lima}

\newcommand{\azero}{\textcolor{blue}{0}}
\newcommand{\aum}{\textcolor{blue}{1}}
\newcommand{\adois}{\textcolor{blue}{2}}
\newcommand{\rzero}{\textcolor{red}{0}}
\newcommand{\rum}{\textcolor{red}{1}}
\newcommand{\rdois}{\textcolor{red}{2}}

\newcommand{\keyw}[1]{{\bf #1}}

\begin{document}


\begin{frontmatter}

	\title{Serial and Unserial Combinatorial Families}

	\author{\melo}
	
	\ead{glaucio.melo@gmail.com}
	
	\ead[url]{http://www.dei.unicap.br/{$\thicksim$}glaucio}
	
	\address{Departamento de Estatística e Informática, 
	Universidade Católica de Pernambuco, Recife, Pernambuco, Brazil}

	\author{\oliveiralima}
	
	\ead{eal@dei.unicap.br}
	
	\ead[url]{http://www.dei.unicap.br/{$\thicksim$}eal}
	
	\address{Departamento de Estatística e Informática, 
	Universidade Católica de Pernambuco, Recife, Pernambuco, Brazil}


	\begin{abstract}
	This article presents the \emph{Serial and Unserial Methods} (SUM). 
	The algorithms are strongly related to the first
	part of a classical reference in combinatorics, the
	\emph{Combinatorial algorithms for computers and calculators},
	from Albert Nijenhuis and Herbert Wilf. The \emph{Serial Method}
	proposal is to obtain the output of a specific kind of
	combinatorial family from its position on the list of all
	combinatorial possibilities. The \emph{Unserial Method} is the
	inverted step of Serial Method, getting the serial number from the
	combinatorial family given as input. The \emph{serial number} is
	the position of the combinatorial family on the
	list.
	\end{abstract}
	
\end{frontmatter}


\section{Introduction}
On the classical reference of combinatorial algorithms \cite{wi},
A. Nijenhuis and H. Wilf presents the \emph{Next} algorithms, and
the proposal of this methods is to obtain the next output of a
combinatorial family from the actual one. Getting an example from
Nijenhuis-Wilf's reference, the \emph{Next Permutation} \cite{wi}
algorithm generates de whole list of permutations without
considering the list position of each permutation contained in
the list. This is a fast way to obtain all the permutations or a
set of permutations from a specific one. This work is a set of algorithms 
which are concentrated to get a specific
combinatorial family on the list. The Serial and Unserial
combinatorial family in this article is the permutation,
composition, partition of an {$n$}-set, {$k$}-subsets and subsets. The combinatorial 
family exposed in this article is divided on the sections below.


\section{Serial Permutation Method (SPM)}


\subsection{Basic Concepts}

The permutation algorithms are classified in two main groups: 
one that creates a set of permutations from the identity
permutation and other that produces a set of permutations by
means of simple changes between the vector elements, creating a
new permutation from the actual one. The algorithm \emph{Next
Permutation for N letters} \cite{wi} belongs to the second group,
creating a complete set of permutations through successive
algorithm invocations, getting the next permutation until reaches
the last vector of the list. The algorithm \emph{Next Permutation}
is creates the next permutation vector using only local
information. The next vector is different from anyone else that 
was defined before.

There is a specific problem on this creating process:
\begin{quote}
\emph{Is it possible to get directly a permutation vector located
on a specific position on the list of vectors, excluding the
alternative of the Next Permutation algorithm's successive
execution until gets the desired vector?}
\end{quote}

The \emph{Serial Permutation Method} (SPM) answers this question,
being able to process the same list of permutation that the
\emph{Next Permutation} algorithm does, with a difference: the SPM
only needs the serial number and the permutation vector's size to
process the desired output. It is also possible to invert this
process from the SPM, getting the serial number through the
permutation vector.


\subsection{SPM's Construction}

The SPM was developed from the observation of the \emph{Next
Permutation's} data outputs. This algorithm has an auxiliary
variable called offset vector (also called inversion vector) \cite{wi}, which is defined below.

Consider {$p$} the permutation vector with {$n$} elements and let
{$d$} be the offset vector with {$n - 1$} elements, the {$d$}
vector is defined by \cite{wi}:

\begin{equation} \label{eq}
d_i = |j: j \leq i, ~p_j > p_{i+1}| ~~i = 1, 2, ... ,n
\end{equation}

That is explained:
\begin{quote}
\emph{How many numbers are there biggest than {$p_{i+1}$} between
the beginning of the permutation vector and the {$i$} index of the
{$d$} vector?}
\end{quote}

Being {$n!$} the complete list of the created vectors, we have in
Table \ref{tab_1} an example of output from the {$p$} permutation
vector related to its serial, for {$n = 5$} and the {$serial =
32$}.

\begin{table}[!hbp]
  \caption{Input and output processing that are 
  common in the \emph{Next Permutation} algorithm.}
  \label{tab_1}
\begin{tabular}{c|c}
 \hline
  Input & Output \\
  \hline
  {$n = 5$} & {$p = (3, 5, 1, 2, 4)$} \\
  {$Serial = 32$} & {$d = (0, 2, 2, 1)$} \\
  \hline
 \end{tabular}
\end{table}

The SPM is subdivided in two steps, and the data input of the
second step corresponds to the data output of the first step:

\begin{enumerate}

\item  Given a serial number and the size of the permutation
vector, determine the offset vector; 

\item Given the offset vector, determine the permutation vector.

\end{enumerate}


\subsection{Algorithm to Attain the Offset Vector}

To solve the first SPM step it is necessary an algorithm to get
the offset vector, which is called \emph{Serial Offset Algorithm}
(SOA). The algorithm does a calculus that reflects the pattern of
the offset vectors' creation throughout the whole list of
permutation, using only local information. In Table \ref{table},
we can visualize the pattern of the offset vector in the whole
list of permutation, where the detached column will serve as a
guide to the determination that the SOA decodes.

\begin{table}
  \caption{Full list of permutations, with {$n = 4$}}
  \label{table}

\begin{tabular}{c|c|c}
 \hline
   Serial  &p  &d\\
 \hline
    1   &(1,2,3,4)  &(0,\azero,0)\\
  \hline
    2   &(2,1,3,4)  &(1,\azero,0)\\
  \hline
    3   &(3,1,2,4)  &(1,\aum,0)\\
  \hline
    4   &(1,3,2,4)  &(0,\aum,0)\\
  \hline
    5   &(2,3,1,4)  &(0,\adois,0)\\
  \hline
    6   &(3,2,1,4)  &(1,\adois,0)\\
  \hline
    7   &(4,2,1,3)  &(1,\rdois,1)\\
  \hline
    8   &(2,4,1,3)  &(0,\rdois,1)\\
  \hline
    9   &(1,4,2,3)  &(0,\rum,1)\\
  \hline
    10  &(4,1,2,3)  &(1,\rum,1)\\
  \hline
    11  &(2,1,4,3)  &(1,\rzero,1)\\
  \hline
    12  &(1,2,4,3)  &(0,\rzero,1)\\
  \hline
    13  &(1,3,4,2)  &(0,\azero,2)\\
  \hline
    14  &(3,1,4,2)  &(1,\azero,2)\\
  \hline
    15  &(4,1,3,2)  &(1,\aum,2)\\
  \hline
    16  &(1,4,3,2)  &(0,\aum,2)\\
  \hline
    17  &(3,4,1,2)  &(0,\adois,2)\\
  \hline
    18  &(4,3,1,2)  &(1,\adois,2)\\
  \hline
    19  &(4,3,2,1)  &(1,\rdois,3)\\
  \hline
    20  &(3,4,2,1)  &(0,\rdois,3)\\
  \hline
    21  &(2,4,3,1)  &(0,\rum,3)\\
  \hline
    22  &(4,2,3,1)  &(1,\rum,3)\\
  \hline
    23  &(3,2,4,1)  &(1,\rzero,3)\\
  \hline
    24  &(2,3,4,1)  &(0,\rzero,3)\\
  \hline
 \end{tabular}
\end{table}

Let {$d$} be the offset vector with {$n$} elements. Each {$d$}
vector's column showed in Table \ref{table} has the following
properties:
\begin{itemize}

\item Let {$k$} be the value of the current column ({$k = 1, 2,
..., n$}), the elements of the {$d_k$} column ({$0 \leq d_k \leq
k$}) repeat {$k!$} times until the element of the column gets the
{$k$} value. This can be visualized on the first part of the
detached column in Table \ref{table}, illustrated in blue, which
the number of detached elements in blue or red is defined by {$(k
+ 1)!$} 

\item When {$d_k = k$} in {$(k+1)!$} times, the following
elements from the list are put on the inverted form, like showed
in red in Table \ref{table}. The direct and inverted list
intercalate themselves until complete all the positions of the
column.
\end{itemize}

Exists throughout each offset vector column an intercalation
between the two kinds of lists (direct and inverted), we can
consider this fact as an element of parity in the list,
considering as an even the list in its direct form and as odd the
list in its inverted form. To determinate the relation between the
serial number and the parity of the list, follows:

\begin{equation}  \label{eqqf}
f = \displaystyle\biggl\lfloor\frac{s-1}{(k+1)!}\displaystyle\biggl\rfloor
~mod~2
\end{equation}

Where:
\begin{itemize}
    \item {$f$}: Determines if the list is direct {$(f = 0)$} or inverted {$(f = 1)$};
    \item {$s$}: Serial number;
    \item {$k$}: Index of the offset vector;
    \item {$\lfloor x \rfloor$}: Floor function. It returns the biggest integer value smaller than {$x$};
    \item {$mod$}: An operation that returns a division's rest.
\end{itemize}

If the list had been direct, we attribute to the elements of the
offset vector:

\begin{equation}    \label{eqt1}
d_k = \displaystyle\biggl\lfloor\frac{[(s-1) ~mod~
(k+1)!]}{k!}\displaystyle\biggl\rfloor
\end{equation}

If the list had been inverted, we attribute the complement which
would be the direct list:

\begin{equation}    \label{eqq}
d_k = k - \displaystyle\biggl\lfloor\frac{[(s-1) ~mod~
(k+1)!]}{k!}\displaystyle\biggl\rfloor
\end{equation}

Each attribution is made through a loop that go through the offset
vector. The equations \ref{eqt1} and \ref{eqq} can be joined in a
unique form, but is a inefficient way to compute the equation:

\begin{equation}    \label{eqj}
d_k = f.k + (-1)^f. \displaystyle\biggl\lfloor\frac{[(s-1) ~mod~
(k+1)!]}{k!}\displaystyle\biggl\rfloor
\end{equation}

Where {$f$} is defined on equation \ref{eqqf}.
\\\\
\textbf{Serial Offset Algorithm (SOA)}
\newline
Routine Specifications:
\begin{itemize}
    \item {$n$}: Size of the offset vector;
    \item {$i$}: Index of the offset vector;
    \item {$d$}: Offset vector, alternating its indices from {$0..n-1$};
    \item {$s$}: Permutation's serial.
\end{itemize}
Routine:
\begin{tabbing}
\quad \=\quad \=\quad \=\quad \=\quad \=\quad \=\quad\=\quad \kill
\keyw{for} {$i \leftarrow 1$} \keyw{to} {$n$} \keyw{do} \\
\> \> \keyw{if} {$ \lfloor (s - 1) ~/~ (i+1)! \rfloor ~mod~ 2 = 1$} \keyw{then}\\
\> \> \> {$d_{i-1} \leftarrow i - \lfloor((s-1) ~mod~ (i+1)!) ~/~i!)\rfloor$}\\
\> \> \keyw{else}\\
\> \> \> {$d_{i-1} \leftarrow \lfloor((s-1) ~mod~ (i+1)!)~/~ i!)\rfloor$}\\
end \keyw{for}\\
\keyw{return} {$d$}.
\end{tabbing}


\subsection{Algorithm to attain the permutation vector from the
offset vector}

After the SOA is computed, the SPM is concluded with the
\emph{Permutation Algorithm by Offset} (PAO). The PAO does the
SPM's second step, returning the desired output. A relevant topic
for the PAO construction is to find the decoding process of the
permutation vector, with only the offset vector being the input.
We know that the offset vector maps the elements of the
permutation vector. In the decoding, we have:

Let {$p$} be the permutation vector and {$d$} the offset vector:
\begin{equation}
    p = (p_1, p_2, ..., p_n) ~~~d = (d_1, d_2, ..., d_{n-1})
\end{equation}

And {$d$} already has its values determining by the SOA. This is
how the decoding process is made: we know that {$d_1$} has its
values discretely included between {$0$} and {$1$}. Obviously,
{$0$} and {$1$} are the unique possible elements for {$d_1$}. From
this information, we can conclude:

\begin{displaymath}
d_1 = \left\{ \begin{array}{ll}
0, & \textrm{if $p_1 < p_2$}\\
1, & \textrm{if $p_2 < p_1$}\\
\end{array} \right.
\end{displaymath}

We do not consider the possibility of equality between the
elements of the permutation vector, because we know that there is
no repeated elements on the vector, and if we organize it, the
difference between them will be only one unit. Like {$d_1$}, the
element {$d_2$} has its values included between {$0$} and {$2$}.
In this and in other cases, we analyze the current inequality from
the inequality that was created previously. For {$d_2$}, we have:

\begin{displaymath}
d_1 = \left\{ \begin{array}{ll} 0, & d_2 = \left\{
\begin{array}{ll}
0, & \textrm{if $p_1 < p_2 < p_3$}\\
1, & \textrm{if $p_1 < p_3 < p_2$}\\
2, & \textrm{if $p_3 < p_1 < p_2$} \\
\end{array} \right.
\\\\
1, & d_2 = \left\{
\begin{array}{ll}
0, & \textrm{if $p_2 < p_1 < p_3$}\\
1, & \textrm{if $p_2 < p_3 < p_1$}\\
2, & \textrm{if $p_3 < p_2 < p_1$} \\
\end{array} \right.\\
\end{array} \right.
\end{displaymath}

We can conclude that as the value of {$d_2$} increases, {$p_3$}
"slides" on the left through the inequality. In general, we have:
\begin{equation}
d_k \in \{0, 1, ..., k\} ~~~ \{p_i < p_j < ... < p_k < ... <
p_t\} ~~~ i, j, t \neq k
\end{equation}

With {$p_k$} in the inequality with {$d_k$} positions, counted
from right to left, because the order of the elements is
ascendent. Finished the offset vector's raster, we will get a set
of inequalities that informs the order of the permutation
elements. We have, for example:

Given a {$d = (0,2,2)$} offset vector, the inequality for a {$p$}
permutation vector is:
\begin{equation} \label{eq:des}
p_3 < p_4 < p_1 < p_2
\end{equation}

Being the last inequality at \ref{eq:des} the final disposition
between the elements of the vector. After the mapping of the
permutation's elements was done, we can say that each element of
the vector ordered in \ref{eq:des} corresponds to the elements of
the identity permutation. This fact classifies the SPM in the
first group described in the basic concepts of this section. On the
\ref{eq:des} inequality, we have:

\begin{equation}
(p_3 = 1) < (p_4 = 2) < (p_1 = 3) < (p_2 = 4)
\end{equation}

The next step to get the final output is to arrange each one of
the elements ordered on its own positions. As the {$p$} vector is
ordered like {$(p_1, p_2, ..., p_n)$}, we have:

\begin{equation}
 p = (p_1, p_2, p_3, p_4)  ~~~ \Rightarrow ~~~ p = (3, 4, 1, 2)
\end{equation}

So we can get the SPM's final output. In implementation terms, to
compute the input built from the permutation vector (from right to
left) the corresponding indices of the vector were used in
relation to its complement. We will notice at the implementation
that in the insertion moment of the current element, if it
replaces another, the elements of inequality located on the left
will "slide" to the left side, allocating space for the current
element. The "slides to the left" operation is implemented on the
subroutine \emph{push}.
\\\\
\textbf{Permutation Algorithm by Offset (PAO)}

Routine Specifications:
\begin{itemize}
    \item {$n$}: Size of the permutation vector;
    \item {$i,j$}: Indices of the algorithm's vectors;
    \item {$p$}: Permutation vector, alternating its elements on {$0..n-1$};
    \item {$r$}: Vector which will keep the element's position before being ordered.
\end{itemize}
Routine:
\begin{tabbing}
\quad \=\quad \=\quad \=\quad \=\quad \=\quad \=\quad\=\quad \kill
{$r_{n-1} \leftarrow 1$}\\
\keyw{for} {$i \leftarrow 1$} \keyw{to} {$n-1$} \keyw{do} \\
\> \keyw{if} {$r_{n-1-d_{i-1}} = 0$} \keyw{then}\\
\> \> {$r_{n-1-d_{i-1}}$} {$\leftarrow i+1$} \\
\> \keyw{else}\\
\> \> {$r_{n-1-d_{i-1}}$} {$\leftarrow push$}\\
end \keyw{for}\\\\
\keyw{for} {$i \leftarrow 0$} \keyw{to} {$n-1$} \keyw{do} \\
\> {$p_{r_i-1} \leftarrow i+1$}\\
end \keyw{for}\\
\keyw{return} {$p$}.\\\\
\textbf{Subroutine Push} \\
\keyw{for} {$j \leftarrow (n-1) - i$} \keyw{to} {$j < (n - 1) - d_{i-1}$} \keyw{do} \\
\>  {$r_j \leftarrow r_{j+1}$}\\
end \keyw{for}\\
\keyw{return} {$i+1$}.\\
\end{tabbing}


\section{The Inverted Process of SPM (Unserial Method)}

We can describe now the inverted process of the SPM's. We have the
permutation vector as the input, and the desired output is the
correspondent serial number. The ingenuous process to get the
serial number for the permutation vector is the raster of the
permutation list, comparing the vectors one by one, until gets the
equivalent vector computed on the input, being the returned value
the loop's index that does this raster. However, we can find the
serial value inverting the SPM's steps:

\begin{enumerate}
    \item Given the permutation vector, find the offset vector;
    \item Given the offset vector, find the serial number.
\end{enumerate}

Like the SPM, the second process depends on the first, with the
first step data output corresponding to the second step data
input.


\subsection{Algorithm to attain the offset vector from the permutation vector}

We will call this algorithm as the \emph{Offset Algorithm by
Permutation} (OAP). It uses on details the definition of the
offset vector \cite{wi} (see equation \ref{eq}). Two nested loops
add the value of each element of the offset vector.
\newline \newline
\textbf{Offset Algorithm by Permutation (OAP)}
\newline
Routine Specifications:
\begin{itemize}
    \item {$n$}: Size of the permutation vector;
    \item {$i, j$}: Indices of algorithms's vectors;
    \item {$p$}: Permutation vector alternating its elements on {$0..n-1$};
    \item {$d$}: Offset vector alternating its elements on {$0..n-2$}.
\end{itemize}
Routine:

\begin{tabbing}
\quad \=\quad \=\quad \=\quad \=\quad \=\quad \=\quad\=\quad \kill
\keyw{for} {$i \leftarrow 0$} \keyw{to} {$n-2$} \keyw{do} \\
\> \keyw{for} {$j \leftarrow 0$} \keyw{to} {$i$} \keyw{do} \\
\> \> \keyw{if} {$p_j > p_{i+1}$} \keyw{then}\\
\> \> \> {$d_i \leftarrow d_i + 1$}\\
\> end \keyw{for}\\
end \keyw{for}\\
\keyw{return} {$d$}.\\
\end{tabbing}


\subsection{Algorithm to attain the serial number from the offset vector}

Like the OAP, the \emph{Serial Algorithm by Offset} (SAO) also
uses definitions already showed at this section. We know that the
value attributed for the elements of the offset vector, defined on
\ref{eqt1} and \ref{eqq} can be ordered depending on the serial.
Given the variables:
\begin{itemize}
    \item {$d = (d_1,d_2,...,d_k,..., d_n)$}: Offset vector;
    \item {$d_k$}: Offset vector's element;
    \item {$s$}: Permutation's serial number;
    \item {$k$}: Index of the offset vector;
    \item {$q$}: Quotient of the division between {$s-1$} and {$(k+1)!$};
    \item {$\lfloor x \rfloor$}: It returns the biggest integer value smaller than {$x$}.
\end{itemize}
For direct list, we have:

\begin{equation} \label{eqdez}
s-1 = d_k . k! + \lfloor q \rfloor . (k+1)!
\end{equation}

For inverted lists, we have:

\begin{equation} \label{eqonze}
s-1 = (k - d_k). k! + \lfloor q \rfloor . (k+1)!
\end{equation}

An important question for the algorithm implementation is how to
find the variable's {$q$} value. To solve this question we have to
consider that {$d_n$} belongs to a column that has only one list
on a direct disposition (see Table \ref{table}, last {$d$} vector
column). With this information, we conclude that the quotient for
this column corresponds to a value between [0,1[, making:

\begin{equation}
\lfloor q \rfloor . (k+1)!
\end{equation}

Corresponds to zero. With this,  \ref{eqdez} and \ref{eqonze} is
equivalent to:

\begin{equation}  \label{eqdoze}
 s-1 = d_k . k!
\end{equation}

For direct lists, and

\begin{equation}   \label{eqtreze}
s-1 = (k - d_k). k!
\end{equation}

For inverted lists.

The strategy of implementation to get the serial is done by
incremental mode, working the current information being based on
past information, doing a raster on the offset vector from right
to left. This is particularly useful on this case, considering the
next value that will converge to the final serial that belongs to
a group related with the previous serial elements. Each quotient
{$q$} during the raster will corresponds to zero, because we are
considering the division in relation to the gotten serial. Another
relevant question is how to classify if the previous element of
the offset vector it is contained in a direct or inverted list.
Looking at Table \ref{table}, we can make easily an equivalence
table, showed in Table \ref{tabletres}.

\begin{table}
  \caption{Next element status list of the offset vector}
  \label{tabletres}
\begin{tabular}{c|c l c}
\hline
\multicolumn{2}{c|}{Previous Value - Parity} && Next Value - Parity\\
\hline
Element & List && List\\
\hline
  Even & Direct & {$\Longrightarrow$} & Direct\\
  Even & Inverted & {$\Longrightarrow$} & Inverted\\
  Odd & Direct & {$\Longrightarrow$} & Inverted\\
  Odd & Inverted & {$\Longrightarrow$} & Direct\\
\hline
\end{tabular}
\end{table}

Exemplifying the attainment serial process from the offset vector,
we have:

Let {$d = (1,0,3)$} be the offset vector and {$s$} the permutation
serial. As the raster is done from right to left, we work first
with the value {$3$}. As was already said, we have as an initial
input a direct list which {$3$} is contained. For direct lists, we
use an equation showed on \ref{eqdoze}:

\begin{equation}
s \leftarrow 3 . 3!  \Rightarrow 18
\end{equation}

After that, we observer the next element of the offset vector. As
the previous element is in a direct list and it is odd, the
element {$0$} will be in a inverted list. Using \ref{eqtreze}, we
have:

\begin{equation}
s \leftarrow s + (2 - 0) . 2! \Rightarrow 18 + 4 \Rightarrow 22
\end{equation}

Finally, observing the last element of the offset vector, we have
the previous element which is in a inverted list and it is even.
Looking at Table \ref{tabletres}, we evidence that the element
{$1$} is in a inverted list. Thus, we have:

\begin{equation}
s \leftarrow s + (1 - 1) . 1! \Rightarrow 22 + 0 \Rightarrow 22.
\end{equation}

After we have calculated the serial in a incremental mode, we add
{$1$} to the serial, because the equations \ref{eqdoze} and
\ref{eqtreze} depend on {$s-1$}. Thus, the serial value
characterizes it between {$(1, n!)$}. Where {$n$} is the size of
the permutation value given as an input:

\begin{equation}
s \leftarrow s + 1 \Rightarrow 22 + 1 \Rightarrow 23.
\end{equation}

The serial {$23$} is the corresponding serial to the offset vector
{$(1,0,3)$} which corresponds to the permutation vector
{$(3,2,4,1)$}. On the algorithm's implementation showed here, the
increment of one unit to the serial was made in the beginning,
before get into the loop. A boolean variable was specified to
determine if the list which belongs the elements is direct or
inverted. The conditions that determine if the list is direct or
inverted is optimized from four to two conditions. The algorithm
returns a non-negative integer value, corresponding to the serial
number required.
\\\\
\textbf{Serial Algorithm by Offset (SAO)}
\newline
Routine Specifications:
\begin{itemize}
    \item {$s$}: Serial value;
    \item {$n$}: Size of the offset vector;
    \item {$i$}: Index of the offset vector;
    \item {$d$}: Offset vector, alternating its elements on {$0..n-1$};
    \item {$direct$}: Boolean variable. It determines if the list is direct or not.
\end{itemize}

\begin{tabbing}
\quad \=\quad \=\quad \=\quad \=\quad \=\quad \=\quad\=\quad \kill
{$s \leftarrow [~d_{n -1} . (n-1)!~] + 1$}\\
{$direct \leftarrow true$}\\
\keyw{for} {$i \leftarrow n-1$} \keyw{to} {$1$} \keyw{do} \\
\> \keyw{if} {$(d_i= even ~and~ direct = true)$} {$or$} {$(d_i = odd ~and~ direct = false)$} \keyw{then}\\
\> \> {$direct \leftarrow true$}\\
\> \> {$s \leftarrow s +  d_{i-1} . (i-1)!$}\\
\> \keyw{else}\\
\> \keyw{if} {$(d_i = odd ~and~ direct = true)$} {$or$} {$(d_i = even ~and~ direct = false)$} \keyw{then}\\
\> \> {$direct \leftarrow false$}\\ 
\> \> {$s \leftarrow s +  (i - d_{i-1}) . (i-1)!$}\\
end \keyw{for}\\
\keyw{return} {$s$}.\\
\end{tabbing}


\section{Serial Composition Method}


\subsection{Basic Concepts} 

On the combinatorial family, the composition of an integer 
{$n$} in {$k$} parts is defined by:

\begin{equation}
n = r_1 + r_2 + \ldots + r_k ~~~ r_i \geq 0 ~~~ i=1..k
\end{equation}

Where the order of the elements is important on the compositions
generations. The \emph{Next Composition Algorithm} \cite{wi} does
the task of making composition, obtaining the next composition
starting from the one before, interactively working to reach until
the last composition of the list. The proposal of the described
method on this section is to obtain the vector from the list
position, with no need of processing the compositions one by one
until getting the expecting vector.


\subsection{Construction of the SCM}

The SCM was built from the repetition patterns which are present
on the composition vector throughout its list. It is possible to
see that this repetition is done in specific positions on the
compositions lists: they can be obtained through calculus. This
feature of the composition vector makes possible the cast among
the specific positions on the list. This way, it decreases the
number of needing interactions to find the vector of a specific
position. The calculus of the specific positions is defined below.


\subsection{Used definitions on the built of SCM}

It is known that the total number of compositions of {$n$} in
{$k$} parts is defined by \cite{wi}:

\begin{equation} \label{eqc1}
J(n,k) = \left(
\begin{array}{ccc}
n + k -1 \\
n\\
\end{array} \right)
\end{equation}
It is also know that:
\begin{equation} \label{eqc2}
\sum_{i=0}^{k} \left(
\begin{array}{ccc}
n - i \\
k - i \\
\end{array} \right) =
\left(
\begin{array}{ccc}
n + 1 \\
k\\
\end{array} \right) \;
\end{equation}

Associating \ref{eqc1} and \ref{eqc2} we have:

\begin{equation} \label{eqc3}
\sum_{i=0}^{n} \left(
\begin{array}{ccc}
n + k - 2 - i \\
n - i \\
\end{array} \right) =
\left(
\begin{array}{ccc}
n + k -1 \\
n\\
\end{array} \right) \;
\end{equation}

For a {$c$} composition vector of a {$k$} length from the {$L$}
composition list, the modification of {$kth$} element in {$L$}
is determined by the index of the partial sum defined on the
equation \ref{eqc3}. This way, it is possible to know when a
component of {$c$} stops to repeat its current element to get
modifications. Taking as a basis equation \ref{eqc3}, we can see
this index modification with the definition of the {$M$} matrix,
exposed in \ref{eqc4}.

\begin{equation} \label{eqc4}
M(n,k) = \left(
\begin{array}{cccc}
\left(\begin{array}{ccc}n + k -2 \\n\\\end{array} \right) &
\left(\begin{array}{ccc}n + k -3 \\n-1\\\end{array} \right) &
\ldots &
\left(\begin{array}{ccc} k - 2 \\0\\\end{array} \right)\\\\
\left(\begin{array}{ccc}n + k -3 \\n\\\end{array} \right) &
\left(\begin{array}{ccc}n + k -4 \\n-1\\\end{array} \right) &
\ldots &
\left(\begin{array}{ccc} k - 3 \\0\\\end{array} \right)\\
\vdots & \vdots & \ddots &\vdots\\
\left(\begin{array}{ccc}n\\n\\\end{array} \right) &
\left(\begin{array}{ccc}n -1\\n-1\\\end{array} \right) & \ldots &
\left(\begin{array}{ccc}0\\0\\\end{array} \right)\\
\end{array} \right)
\end{equation}

\begin{itemize}
\item Each line of {$M$} that corresponds to the elements of the sum on equation \ref{eqc3};

\item From one line to another, the superior element of the binomial is decreased in one unity;

\item {$M$} has {$k-1$} lines and {$n+1$} columns.

\end{itemize}

Now let's see bellow how {$M$} is used to represent the SCM
execution.
\subsection{Description of the {$M$} Matrix Raster}
Initially, we have as initial information the serial number of the
composition on the list, to obtain the composition vector. Let's
suppose {$s$} is the serial related to the input variable of the
method to a composition vector with a {$k$} length. The elements
of the first line of {$M$} are added until this sum is over
{$s-1$}. When it's found a value in the sum which is over {$s-1$},
the corresponding value of the column which the added element
immediately before the actual one is found is given to the last
element of the composition vector found. After finding the last
element (the attributions of each element of the composition
vector are done from the end to the beginning until the second
element of the composition vector), the next element, that means
the element before the last one of the composition vector is found
in {$M$} going down in diagonal to the next line of {$M$},
re-starting the counting of indices starting from the actual
column. This is done in a successive way, until the sum of the
elements which were visited in {$M$} is {$s-1$}. Let's see an
example of this strategy in a matrix {$M(7,5)$} with serial {$s-1
= 282$}, seen in \ref{eqc5}. The result of the raster is the
composition vector {$c(283)= \{1,0,2,1,3\}$}. It is known that the
first element of the vector does not need to be calculated using
{$M$}, once it can be obtained from the complement of the sum of
the elements already found through the raster done in {$M$}. There
is in \ref{eqc5} on the first row with the partial sum until the
third element (we allocate the number {$3$} on the last position
of the list), the second row with only one element of the partial
sum (we allocate the number {$1$} on the position before the last
one of the list), and so on. The sum of the underlined numbers in
\ref{eqc5} converged to {$s - 1$}.

\begin{equation} \label{eqc5}
	M(7,5) = \left(
	\begin{array}{cccccccc}
	\textbf{\underline{120}} & 
	\textbf{\underline{84}} & 
	\textbf{\underline{56}} & 35 & 20 & 10 & 4 & 1 \\
	36 & 28 & 21 & \textbf{\underline{15}} & 10 & 6 & 3 & 1\\
	8 & 7 & 6 & 5 & \textbf{\underline{4}} & \textbf{\underline{3}} & 2 & 1 \\
	1 & 1 & 1 & 1 & 1 & 1 & 1 & 1 \\
\end{array} \right)
\end{equation}

\subsection{Description of the SCM}

Related to implementation, the algorithm proposed here to the SCM
is abstract to the structure of {$M$} matrix, defined above. The
loop is traced as defined, applying on the first element of the
composition vector the corresponding value of the complement of
the sum of the elements yet to come. On the subroutine
\emph{element} there is the procedure to obtain each element of
the composition vector.\\
\\
\textbf{Serial Composition Method (SCM)}
\\
Algorithm specifications:
\begin{itemize}
\item {$n$}: Number of the composition (Composition of {$n$} elements in {$k$} parts);
\item {$k$}: Parts of the composition;
\item {$s$}: Serial of the composition;
\item {$a$}: Auxiliary variable which makes the increment
to the convergence of the serial number;
\item {$x,y$}: Auxiliary variables, for the definition of new binomial indices;
\item {$C_{i,j}$}: Combination of {$i$} elements {$j$} by {$j$};
\item {$z$}: Value of the complement, used to apply the value to the 
first element of the composition vector.
\end{itemize}
Routine:

\begin{tabbing}
\quad \=\quad \=\quad \=\quad \=\quad \=\quad \=\quad\=\quad \kill
{$z \leftarrow a \leftarrow 0$} \\
{$x \leftarrow n + k - 2$}\\
{$y \leftarrow n$}\\
\keyw{for} {$i \leftarrow 0$} \keyw{to} {$k-2$} \keyw{do} \\
\> {$c_{k-1-i} \leftarrow element $}\\
\> {$z \leftarrow z + c_{k-1-i}$}\\
end \keyw{for}\\
{$c_0 \leftarrow n - z$}\\
\keyw{return} {$c$}.\\
\\
Subroutine element \\
\keyw{for} {$j \leftarrow 0$} \keyw{to} {$n-1$} \keyw{do} \\
\> \keyw{if} {$a + C_{x-i-j,y-j} \leq s-1$} \keyw{then}\\
\> \> {$a \leftarrow a +  C_{x-i-j,y-j}$}\\
\> \keyw{else}\\
\> \> {$x \leftarrow x - j$}\\
\> \> {$y \leftarrow y - j$}\\
\> \> \keyw{return} {$j$}.\\
end \keyw{for}\\
\keyw{return} {$n$}.\\
\end{tabbing}


\section{The inverse process of SCM (Unserial Method)}

Related to the inverse process of the SCM, the actual view is the
composition serial number obtaining, having as input data the
composition vector of a number {$n$} in {$k$} parts. In this case,
each component of the composition vector is seen as part of
superior intervals of a nested loops, which makes the raster on
the composition vector. The inverse process has concise data input
related to the data input of the SCM, once the inverse process
will not work under the value convergence of a serial number
method: everything is already well defined on the composition
vector components, known that is only necessary to process the
loops with related interactions to each component of the
composition vector give as input. As to the SCM, the proposed
algorithm for the inverse process to the SCM works only with
indices referred to the defined matrix on SCM. The use of the
indices abstracting from the matrix structure decreases the memory
use and processes only what is needed for the calculus of the
trace to be run on {$M$} matrix.\\
\\
\textbf{Serial Composition Method (Inverse process)}
\\
Algorithm specifications:
\begin{itemize}
\item {$n$}: Number of the composition (Composition of {$n$} elements in {$k$} parts);
\item {$k$}: Parts of the composition;
\item {$s$}: Serial of the composition;
\item {$x,y$}: Auxiliary variables, for the definition of new binomial indices;
\item {$C_{i,j}$}: Combination of {$i$} elements {$j$} by {$j$}.
\end{itemize}

Routine:

\begin{tabbing}
\quad \=\quad \=\quad \=\quad \=\quad \=\quad \=\quad\=\quad \kill
{$x \leftarrow n + k -2$} \\
{$y \leftarrow n$} \\
{$s \leftarrow 1$} \\
\keyw{for} {$i \leftarrow k-1$} \keyw{to} {$1$} \keyw{do} \\
\> \keyw{for}  {$ j \leftarrow 0$} \keyw{to} {$c_i -1$} \keyw{do} \\
\> \> {$s \leftarrow s + C_{x-j-[(k-1)-i],y-j}$} \\
\> end \keyw{for}\\
\> {$x \leftarrow x - c_i$} \\
\> {$y \leftarrow y - c_i$} \\
end \keyw{for}\\
\keyw{return} {$s$}.\\
\end{tabbing}


\section{Serial Partition of an n-Set Method}


\subsection{Basic Concepts}

For the sets partitions, we consider a family of subsets {$T_1,
T_2,\ldots,T_k $} contained in a set {$S = \{1,2,\ldots,n\}$} that
satisfy the conditions \cite{wi}:

\begin{equation}
T_i \cap T_j = \oslash ~~~ (i \neq j)
\end{equation}

\begin{equation}
\bigcup _{i=1}^{k} T_i = S
\end{equation}

\begin{equation}
T_i \neq \oslash \> (i = 1, 2, \ldots,k)
\end{equation}

We don't consider the order of the elements contained on each one
of {$S$} subsets. The algorithm \emph{Next Partition of an
n-Set} \cite{wi} finds the next partition set from the current one,
working only with local information. The section's proposal is to
obtain a partition of {$S$} related to its position on the list of
partition, without considering the information about the partition
vector.

\subsection{SPSM's Construction}

The SPSM was built from the analysis of the data output of the
\emph{Next Partition of an {$n$}-Set} algorithm. Each index that
denotes the owned elements to a subset has a pattern model and the
search for these indices is made by a combinatorial structure that
take as basis a tree that represents all the partitions of a set
with {$n$} elements, called Bell Tree. The number of nodes of this
structure grows quickly when {$n$} increases. So, the solution for
such problem is to define a new structure that does a mapping of
the tree specific positions. As a reference we took a matrix
structure to keep those positions, called {$D$} Matrix. The Bell
Tree and {$D$} Matrix are defined below.


\subsection{Bell Trees} 

We can define a Bell Number as a number of
partitions possibilities of a set with {$n$} elements. Such number
is defined by \cite{gkp}:

\begin{equation} \label{eqp1}
B_n = \sum_{k=1}^{n} \bigg\{{n\atop k}\bigg\}
\end{equation}

Where {$\big\{{n\atop k}\big\}$} is the Stirling number of second
kind, which is defined by \cite{gkp}:

\begin{equation} \label{eqp2}
\bigg\{{n\atop k}\bigg\} = \frac{1}{k!} \sum_{i=0}^{k-1}(-1)^i {k
\choose i}(k-i)^n
\end{equation}

With the definitions \ref{eqp1} and \ref{eqp2} above, this
subsection proposes to show a combinatorial structure which is
called here \emph{Bell Tree} and the use of this structure to
solve the \emph{Partition of an n-Set} problem. It was used as
basis the definition of a partition tree, shown in the algorithm
\emph{Next Partition of an n-Set} \cite{wi}, associating the tree
nodes with the equations which its terms are defined on the tree
construction properties.


\subsection{Properties of Bell Trees}

Let {$n$} be the number of elements of a set that will be
partitioned and {$B$} the Bell Numbers that denotes a set of terms
that can form a mathematical expression in a node of the tree, the
structure here defined adopt the partition tree structure with
{$n$}-Sets and added to it others properties. Follows below some
properties of the structure:

\begin{itemize}

\item The Tree has {$n$} levels;

\item The Tree nodes are non-negative integers numbers which is
formed by expressions that evolves exclusively Bell Numbers and
multiplicative integer constants;

\item Let {$S$} be a subset contained in {$B$} which denotes the
number of distinct terms of an expression with a node {$N$}, we
determine the number of descendants produced by {$N$} from the
number of elements of {$S$} increased one unit;

\item The number {$R$} of the nodes on a Bell Tree corresponds to:

\begin{equation} \label{eqp3}
R = \sum_{j=1}^{n}B_j
\end{equation}

\item Let {$k$} be the number of descendants of the actual node
{$N$}, the value for all its {$k-1$} are defined by:

Let {$E_{w,z}$} be a mathematical expression with a {$w$} index and {$z$} terms
from the ascendant node, the descendant node corresponds to {$E_{w-1,z}$}. To
illustrate such fact, if we have a ascendant node represented by
the expression:

\begin{equation} \label{eqp4}
E_{w,z} = (B_w - B_{w-1}) - 2(B_{w-1} - B_{w-2})
\end{equation}

Then, the descendants nodes will correspond to the expression:

\begin{equation} \label{eqp5}
E_{w-1,z} = (B_{w-1} - B_{w-2}) - 2(B_{w-2} - B_{w-3})
\end{equation}

For the {$kth$} descendant of {$N$}, we have the expression:

\begin{equation} \label{eq6}
E_{k,z+1} = E_{w,z} - (k-1)E_{w-1,z}
\end{equation}

That it is equivalent to say that the {$kth$} descendant is
equal to the ascendant node value minus {$k-1$} times the value of
the others descendant nodes.

\end{itemize}

\begin{figure}
		\centering
	  \includegraphics[width=290pt,height=80pt]{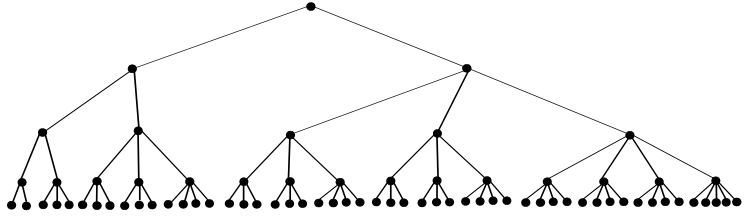}
	  \caption{Bell Tree with level 5}\label{a1}
	\end{figure}

We can now make a tree that follows this formation law. For {$n =
5$}, we have a tree illustrated on Figure \ref{a1}. The three
first levels expressions of this tree are in Table \ref{t1}, which
the number of nodes are entered from top to bottom and from left
to right.

\begin{table}
  \caption{Nodes values of the Tree illustrated on Figure \ref{a1}}
  \label{t1}
\begin{tabular}{l|l}
 \hline
   Nodes  & Expressions  \\
 \hline
    1   & {$B_5$}\\
  \hline
    2   & {$B_4$}\\
  \hline
    3   & {$B_5 - B_4$}\\
  \hline
    4   & {$B_3$}\\
  \hline
    5   & {$B_4 - B_3$}\\
  \hline
    6   & {$B_4 - B_3$}\\
  \hline
    7   & {$B_4 - B_3$}\\
  \hline
    8   & {$(B_5 - B_4) - 2(B_4 - B_3$)}\\
  \hline
 \end{tabular}
\end{table}

\subsection{Definition of {$D$} Matrix}

The {$D$} Matrix is a superior triangular matrix with {$n$} x
{$n$} dimension. Its first column is made by the Bell Numbers

\begin{equation}
    D_{v,0} = B_{n-v}
\end{equation}

and the others columns are defined by the equation below:

\begin{equation}
    D_{u,v} = D_{u,v-1} - v.D_{u+1,v-1}
\end{equation}

The matrix showed on \ref{md} represents the {$D$} Matrix for {$n=6$}.

\begin{equation} \label{md}
D = \left(\begin{array}{c c c c c c}
203 & 151 & 77 & 26 & 6 & 1 \\
52  & 37  & 17 & 5  & 1 &     \\
15  & 10  & 4  & 1  &   &     \\
5   & 3   & 1  &    &     &   \\
2   & 1   &      &      &     &   \\
1   &         &      &      &       &   \\
\end{array}\right)
\end{equation}


\subsection{The {$D$} Matrix fulfilling}

Before process the SPSM, it is necessary fill in the components of
{$D$} Matrix to speed up the attainment calculus of the search on
the Bell Tree. The use of {$D$} Matrix abstracts the construction
of the whole Bell Tree, processing only the necessary to the SPSM
use. If the tree had been built completely, we would have the
number of nodes {$R$} equivalent to the equation \ref{eqp3}.
Whereas on the {$D$} Matrix , we have

\begin{equation} \label{eqpb2}
R = \frac{n^2 + n}{2}
\end{equation}

\textbf{Algorithm for {$D$} Matrix fulfilling}\\
Algorithm Specifications:

\begin{itemize}
\item {$D$}: Matrix that keep the specific positions of the Bell
Tree;

\item {$i,j$}: Indices for the raster of {$D$} Matrix;

\item{$B_n$}: Bell Number.

\end{itemize}
Routine:

\begin{tabbing}
\quad \=\quad \=\quad \=\quad \=\quad \=\quad \=\quad\=\quad \kill
\keyw{for} {$i \leftarrow 0$} \keyw{to} {$n-1$} \keyw{do} \\
\> {$D_{i,0} \leftarrow B_{n-i}$} \\
end \keyw{for}\\
\keyw{for} {$i \leftarrow 1$} \keyw{to} {$n-1$} \keyw{do} \\
\> \keyw{for} {$j \leftarrow 0$} \keyw{to} {$n-i-1$} \keyw{do} \\
\> \> {$D_{j,i} \leftarrow D_{j,i-1} - i.D_{j+1,i-1}$}\\
\> end \keyw{for}\\
end \keyw{for}\\
\keyw{return} {$D$}.
\end{tabbing}


\subsection{SPSM's Specifications}

After that the {$D$} Matrix has its values filled in, the SPSM can
be invoked. We have an external loop that runs all the partition
vector, attributing for each component by means of the sub-routine
\emph{element} that does a raster on the current tree level (the
structure was abstracted from the method). On this raster, we have
a condition that examines if the variable used for the indication
of extrapolation surpassed the serial number token as an input. If
it exceed this value, we go down a tree level; if not, we add the
descendants' current value to the control variable for later
verification of extrapolation, token as an input. If the condition
had not been satisfied for the whole loop, it means that the
search on the tree's level arrived to the last descendant of this
one, indicating that the search will be expanded to the last
descendant of the current tree's level.\\
\\
\textbf{Serial Partition of an {$n$}-Set Method}\\
Algorithm Specifications:

\begin{itemize}

\item {$p$}: Partition vector;

\item {$D$}: Matrix that keep the specific positions of the tree;

\item {$i,j$}: Indices for the raster of {$D$} Matrix;

\item {$k$}: Index for the raster on the {$p$} vector;

\item {$t$}: Index that determines the element of each component
of partition vector;

\item {$a$}: Number that does the convergence for the serial
number given as input;

\item {$s$}: Partition's serial.

\end{itemize}
Routine:
\begin{tabbing}
\quad \=\quad \=\quad \=\quad \=\quad \=\quad \=\quad\=\quad \kill
{$i \leftarrow j \leftarrow a \leftarrow 0$}\\
\keyw{for} {$k \leftarrow 0$} \keyw{to} {$n-1$} \keyw{do} \\
\>    {$p_k \leftarrow element$}\\
end \keyw{for}\\
\keyw{return} {$p$}.\\\\
Sub-Routine {$element$}\\
\keyw{for} {$t \leftarrow 0$} \keyw{to} {$j$} \keyw{do} \\
\> \keyw{if} {$a + D_{i,j} \geq s$} \keyw{then} \\
\> \> {$i \leftarrow i + 1$}\\
\> \> \keyw{return} {$t$}.\\
\> \keyw{else}\\
\> \> {$a \leftarrow a + D_{i,j}$}\\
end \keyw{for}\\
{$j \leftarrow j + 1$} \\
\keyw{return} {$j$}.
\end{tabbing}


\subsection{Algorithm for stylized output of the
partition vectors} 

We know that the output gotten on both SPSM and
Next Partition of an {$n$}-Set corresponds to the indices of each
subset of the partition on the set. For the partition vector on
position {$26$}, we have the output {$(0,1,1,0,0)$} and the
stylized output corresponds to {$(1,4,5)(2,3)$}, indicating that
{$1,4$} and {$5$} is on the first subset (indicated as {$0$}),
with {$2$} and {$3$} on the second subset (indicated as {$1$} on
the partition vector). For show the stylized output, it was used a
\emph{string} vector to get the elements adequately on each subset
which it belongs. Next, the algorithm does the output treatment
and the all elements of each subset be ordered related to the last
elements of partition as an input for the stylized output of the
{$n$}-sets' partitions method.\\
\\
\textbf{Algorithm for stylized output data}\\
Algorithm Specifications:
\begin{itemize}
\item {$s$}: {$Strings$} vector that does the mapping of the
elements referring to the partition vector's indices;

\item {$p$}: Partition vector;

\item {$i$}: Indices for the raster of {$s$} vector;

\item {$r$}: Final output on {$string$} form;

\item {$n$}: Size of the partition vector;

\item {$+$}: Operator that denotes a concatenation between
{$strings$};

\item {$k$}: {$String$} that get the current element of the
stylized output;

\item {$length(w)$}: Function that returns the size of the
{$string$} {$w$};

\item {$substring(w,ini,sup)$}: Function that returns a {$w$}
{$substring$} of the interval between {$ini$} and {$sup$}.

\end{itemize}
Routine:
\begin{tabbing}
\quad \=\quad \=\quad \=\quad \=\quad \=\quad \=\quad\=\quad \kill
\keyw{for} {$i \leftarrow 0$} \keyw{to} {$n-1$} \keyw{do} \\
\>    {$s_{p_i} \leftarrow s_{p_i} + (i + 1) + ","$}\\
end \keyw{for}\\
\keyw{for} {$i \leftarrow 0$} \keyw{to} {$n-1$} \keyw{do} \\
\> \keyw{if} {$length(s_i) > 1$} \keyw{then} \\
\> \> {$k \leftarrow substring(s_i$}, {$0$}, {$length(s_i) - 1)$}\\
\> {$s_i \leftarrow "(" + k + ")"$}\\
\> \keyw{if} {$s_i = "()"$} \keyw{then}\\
\> \> {$r \leftarrow r + s_i$}\\
end \keyw{for}\\
\keyw{return} {$r$}.
\end{tabbing}


\section{The Inverted Process of SPSM (Unserial Method)}

For the inverted process of SPSM, we have a partition vector as
input and the serial number as output. Such process is made taking
each partition vector's component as the number of loops that can
does a sum of each specific position referring to the current node
of the tree mapping through the {$D$} Matrix.


\subsection{Inverted Process of SPSM Specifications}

Initially, we have the indices attribution referring to the {$D$}
Matrix. The {$u$} index has the attribution equals to {$1$}
because it is not necessary do the mapping on the first component
of the {$D$} Matrix, once it always begins with zero for being the
main tree descendant. After that, we have an external loop that
runs the whole partition vector given as input, with another loop
that does the sum to the main positions' serial mapped on {$D$}
Matrix. After the sum attributions to the serial, we have the
checking if the actual component should go down a line on {$D$}
Matrix (equivalent to go down a level on the tree). If it should
go down a level, the line is added by one unit. If not, the column
on {$D$} is added by one unit, indicating the search to the
descendants on the current node reached its last one and the
search be doing on the other tree ramification.\\
\\
\textbf{Inverted Process of SPSM}\\
Algorithm Specifications:
\begin{itemize}
\item {$p$}: Partition vector, given as input;

\item {$i,j$}: Indices referring to the raster on {$p$} vector;

\item {$D$}: Matrix of specific positions on the Bell Tree;

\item {$u,v$}: Indices referring to the mapping on {$D$};

\item {$s$}: Partition's serial, that is the result expected from
de algorithm.
\end{itemize}
Routine:

\begin{tabbing}
\quad \=\quad \=\quad \=\quad \=\quad \=\quad \=\quad\=\quad \kill
{$i \leftarrow j \leftarrow v \leftarrow 0$}\\
{$s \leftarrow u \leftarrow 1$}\\

\keyw{for} {$i \leftarrow 1$} \keyw{to} {$n-1$} \keyw{do} \\
\> \keyw{for} {$j \leftarrow 0$} \keyw{to} {$p_i - 1$} \keyw{do} \\
\> \> {$s \leftarrow s + D_{u,v}$}\\
\>       end \keyw{for}\\
\> \keyw{if} {$j \leq p_{i-1}$} \keyw{then}\\
\> \> {$u \leftarrow u + 1$}\\
\> \keyw{else}\\
\> \> {$v \leftarrow v + 1$}\\
end \keyw{for}\\
\keyw{return} {$s$}.
\end{tabbing}


\section{Serial k-Subset of an n-Set Method}


\subsection{Basic Concepts}

We call {$\big({n \atop k}\big)$} the number of
possibilities to combine {$n$} things on {$k$} different parts. On
analyzed literature \cite{wi}, we have two ways to do this work on
a sequential form. The first method builds the {$k$}-subsets in a
lexicographic order and the second method obtains the next subset
from its predecessor, subtracting one element of the set and
adding on other element of the subset.

The proposal of this section is presents a method that, put in
action iteratively, it creates a list in a lexicographic order
given as input the position on the list of possible combinations.
This strategy is more efficient when we desire to get a
combination on a specific position inside the whole list of
combinations.


\subsection{Getting the combinations on lexicographic order}

The algorithm \emph{Next {$k$}-subset of an {$n$}-Set} \cite{wi} is
able to create on a simple way the combinations in lexicographic
order in a non-recursive way. The recursive model of this
algorithm will be shown on this section optionally. The current
combination of the recursive model is showed through the
{$showOutput$} method, that is indicating here an output model of
generic data for a combination that is going to have as output
elements that can vary between {$0$} and {$n-1$}. The call to
begin the routine must be done on a {$Combine(0)$} way, taking as
stopped of the recurrence, the moment that the parameter did {$k$}
recurrences. Next, we have a recursive model of the algorithm,
without taking longer on its construction.\\
\\
\textbf{Recursive {$k$}-subset of an {$n$}-Set Algorithm}
\\
Algorithm Specifications:
\begin{itemize}

\item {$k$}: Dimension of the subset;

\item {$s$}: Subset vector;

\item {$n$}: Set's cardinality;

\item {$y$}: Vector that determines the changes of {$n$} set
elements in subset {$k$}.

\end{itemize}
Routine:

\begin{tabbing}
\quad \=\quad \=\quad \=\quad \=\quad \=\quad \=\quad\=\quad \kill
{$Combine(i)$}\\
\keyw{for} {$s_i \leftarrow sum(i,i,0)$} \keyw{to} {$n-(k-i)$} \keyw{do} \\
\> \keyw{if} {$i \neq k-1$} \keyw{then}\\
\> \> {$Combine(i+1)$}\\
\> \keyw{else}\\
\> \> {$ShowOutput$}\\
end \keyw{for}\\
{$y_{i+1} \leftarrow 0$}\\
{$y_i \leftarrow y_i + 1$}\\\\
Subroutine {$sum(w,j,z)$}\\
\keyw{for} {$i \leftarrow 0$} \keyw{to} {$j$} \keyw{do} \\
\> {$z \leftarrow z +  y_i$} \\
end \keyw{for} \\
\keyw{return} {$z+w$}.
\end{tabbing}


\subsection{SKSM's Construction}

To build the SKSM, we observe the data output of the subsets' list
created by the \emph{Next {$k$}-subset of an {$n$}-Set} algorithm,
characterizing the repetition's pattern in the elements of subset.
In this case, the pattern can be delineated under a tree model.
The formation law of this tree is showed below.
\subsection{Definition of the Binomial Tree}
The structure that represents the repetition's pattern is
characterized as a tree formed exclusively by binomial
coefficients. Let {$\big({n\atop k}\big)_{w}$} be the current node
of the tree with {$w$} label, its descendants are defined by

\begin{equation}
\bigg({n\atop k}\bigg)_{w} \rightarrow \left\{
\begin{array}{l}
\Big({n-1\atop k-1}\Big)_{w+1}\\
\\
\Big({n-2\atop k-1}\Big)_{w+2}\\
~~~ \vdots \\
\Big({n-k+1\atop k-1}\Big)_{w + n - k + 1}\\
\end{array} \right.
\end{equation}

Each ascendant will have {$n + k -1$} descendants, and the
repetition's pattern of combinations are analyzed through the
insertion of labels in each node of the tree. For each label of
the current node, the descendants nodes enter its labels in
relation to the ascendant node, indicating the value from each
element of the subset founded on the tree.


\subsection{SKSM Specifications}

The SKSM abstracts the construction of the tree, doing the search
calculating only the binomial coefficients and its relations with
the nodes' labels that were visited on the search. The method does
an external loop attributing to each element of the subset the
result of the subroutine {$element$}. The subroutine {$element$}
does a search on the tree's nodes. The abstraction of the tree is
made through the indices changes of the binomial coefficient
attributed on {$x$} and {$y$}. The method is described below.\\
\\
\textbf{Serial {$k$}-Subset of an {$n$}-Set Method}
\\
Algorithm Specifications:
\begin{itemize}

\item {$p$}: {$k$}-dimensional Subset;

\item {$n$}: Cardinality of the Set;

\item {$a$}: Auxiliary variable that is used to check the stopped
of the method;

\item {$x,y$}: Indices of the binomial coefficients functioning as
the formation law of the Binomial Tree;

\item {$r$}: Variable that controls the labels of Binomial Tree;

\item {$s$}: Serial of {$p$} subset;

\item {$C_{n,k}$}: {$\big({n\atop k}\big)$}.

\end{itemize}

Routine:

\begin{tabbing}
\quad \=\quad \=\quad \=\quad \=\quad \=\quad \=\quad\=\quad \kill
{$x \leftarrow n$} \\
{$y \leftarrow k-1$}\\
{$a \leftarrow r \leftarrow 0$}\\
\keyw{for} {$i \leftarrow 0$} \keyw{to} {$k-1$} \keyw{do} \\
\> {$p_i \leftarrow element$}\\
end \keyw{for}\\
\keyw{return} {$p$}.\\\\
Subroutine {$element$}\\
\keyw{for} {$j \leftarrow 1 $} \keyw{to} {$x - y + 1$} \keyw{do} \\
\> \keyw{if} {$a + C_{x-j,y} < s$} \keyw{then}\\
\> \> {$a \leftarrow a + C_{x-j,y}$}\\
\> \keyw{else}\\
\> \> {$x \leftarrow x - j$}\\
\> \> {$y \leftarrow y - 1$}\\
\> \> {$r \leftarrow r + j$}\\
\> \> \keyw{return} {$r$}.\\
end \keyw{for}\\
\keyw{return} {$r$}.\\
\end{tabbing}


\section{The Inverted process of SKSM (Unserial Method)}

As many others serial methods from others combinatorial problems
with same importance, the SKSM also has its inverted process. What
is giving as input is the vector that represents the
{$k$}-dimensional subset with the cardinality of the set, called
{$n$}, getting as output its position (i.e. the serial number) on
the subsets list in a lexicographic order.


\subsection{SKSM Inverse Process Specifications}
On the inverted SKSM process, it is taken as the maximum limit of
the internal loop is the difference between the elements of the
subset. This will determine how much the loop will repeat related
to the elements of the subset. The indices of the binomial
coefficient do an appropriate control for the adding among the
nodes be made correctly. The inverted process also abstract the
tree's structure, getting the result by means of local information
related to the subset given as input.\\
\\
\textbf{Serial {$k$}-Subset of an {$n$}-Set Method (Inverted
Process)}
\\
Algorithm Specifications:
\begin{itemize}
\item {$p$}: {$k$}-dimensional Subset;

\item {$n$}: Cardinality of the Set;

\item {$x,y$}: Indices of the Binomial Coefficients;

\item {$r$}: Variable that will control the internal loop of the
method;

\item {$s$}: Serial of {$p$} subset;

\item {$C_{n,k}$}: {$\big({n\atop k}\big)$}.

\end{itemize}
Routine:

\begin{tabbing}
\quad \=\quad \=\quad \=\quad \=\quad \=\quad \=\quad\=\quad \kill
{$x \leftarrow n$}\\
{$y \leftarrow k-1$}\\
{$s \leftarrow 1$}\\
{$r \leftarrow 0$}\\
\keyw{for} {$i \leftarrow 0$} \keyw{to} {$k-1$} \keyw{do} \\
\> \keyw{for} {$j \leftarrow 1$} \keyw{to} {$p_i - r - 1$} \keyw{do} \\
\> \> {$s \leftarrow s + C_{x-j,y}$}\\
\> end \keyw{for}\\
\>  {$x \leftarrow x - (p_i - r)$}\\
\>  {$y \leftarrow y - 1$}\\
\>  {$r \leftarrow p_i$}\\
end \keyw{for}\\
\keyw{return} {$s$}.
\end{tabbing}


\section{Serial Subset Method}


\subsection{Basic Concepts}

We have {$2^n$} possible configurations to get the subsets of a
{$\{1,2,\ldots,n\}$} set. The disposition of the subset's elements
can be represented by a flag which activates its insertion into
the subset. For example, the configuration {$\{1,0,1,1,0\}$}
represents the subset {$\{1,3,4\}$}. The \emph{Next Subset of an
{$n$}-Set} \cite{wi} algorithm does this job in a sequential form,
returning the next subset from the actual one.

This section proposes to present a method that returns a subset of
an {$n$}-Set from its serial number. The inverted process for this
method is also showed, getting the serial number from the subset
given as input.


\subsection{SSM's Construction} 

The SSM was built from the analysis of repetition's pattern that 
the subsets presents in the complete list of subsets. 
Similarly to the \emph{Serial
Permutation Method} related to repetition's pattern of the offset
vector, the SSM has a regular repetition pattern, being able to
get each subset's component through a closed equation.

Let {$p$} be a subset of a {$n$}-Set, each subset's component with
index {$k = 0,1,\ldots,n-1$} related to a serial number {$s =
1,2,\ldots,2^n$} is defined by

\begin{equation}    \label{eqst}
p_k = \displaystyle\biggl\lfloor\frac{(s-1 + 2^k) ~mod~
2^{k+2}}{2^{k+1}}\displaystyle\biggl\rfloor
\end{equation}

The equation \ref{eqst} doesn't present restrictions in the
pattern's representation. Each component {$p_k \in \{0,1\}$}, and
{$2^k$} in {$(s-1 + 2^k)$} is related with a cyclic structure
which the numbers are presented at the list. This peculiarity is
easily identified on output data of \emph{Next Subset of an
{$n$}-Set} \cite{wi} algorithm and the method here proposed.


\subsection{SSM Specifications}

About the SSM's specifications, we have a loop that do an
association between the equation \ref{eqst} and the structure that
represents the subset. The method returns a subset from the serial
number given as input, where the number {$1$} on the output data
indicates which subset's component will be activated, as explained
in the basic concepts on this section.\\
\\
\textbf{Serial Subset Method}
\\
Algorithm Specifications:
\begin{itemize}

\item {$p$}: Subset of {$n$}-Set;

\item {$n$}: Cardinality of the Set;

\item {$s$}: Serial of {$p$} subset;

\end{itemize}
Routine:

\begin{tabbing}
\quad \=\quad \=\quad \=\quad \=\quad \=\quad \=\quad\=\quad \kill
\keyw{for} {$i \leftarrow 0$} \keyw{to} {$n-1$} \keyw{do} \\
\>  {$p_i \leftarrow \lfloor ((s-1 + 2^i) ~mod~
2^{i+2})~/~2^{i+1}\rfloor$}\\
end \keyw{for}\\
\keyw{return} {$p$}.
\end{tabbing}


\section{The Inverted Method of SSM (Unserial Method)}

In this section the inverse process of the SSM will be shown. The
input data is the subset (the elements of the set that will be
activated) and the output data is the corresponding position in
the list of subsets. In the case at hands, we determine which it
is the moment to walk through specific positions of the identified
repetition's pattern in subsets. The inverted process for subsets
also have similarities with the inverted process of \emph{Serial
Permutation Method}. The list of subsets, from right to left,
presents a sequence of {$0$}'s and {$1$}'s, in this order. From
there, the list for analysis of the subset's next component can be
inverted in the case of this subset's element be equals to {$1$}.
In this case, we cannot advance in the list for convergence of the
serial's desired value. Otherwise, the serial value is modified on
iterative form among the specific position of the list until reach
the subset's serial.


\subsection{SSM Unserial Method Specifications}

In the SSM's specification, it was defined a logical variable
which defines if the list of current subset's component is direct
or inverted, defining as {$true$} for direct lists and {$false$}
for inverted lists, starting this variable as {$true$} because the
components' disposition has been analysed from right to left. The
initial serial's value is {$1$}. After that, we have a decreased
loop that does the verifications in each component of the subset
for determines if it advances or not for the convergence of the
final result. It was observed that such condition in the way as
was constructed can be represented by \emph{eXclusive OR} logical
connective, usually denoted by {$xor$} operator. When the {$xor$}
is satisfied, the logical variable is activated as {$true$},
indicating that the next component is included on a direct list.
Otherwise, the serial advances other specific position
and the list is configured as inverted.\\
\\
\textbf{Inverted Process of the SSM (Unserial Method)}
\\
Algorithm Specifications:
\begin{itemize}

\item {$p$}: Subset of {$n$}-Set;

\item {$n$}: Cardinality of the Set;

\item {$s$}: Serial of {$p$} subset;

\item {$d$}: Logical variable, which indicates if the list of
subsets is on direct or inverted order.

\end{itemize}
Routine:

\begin{tabbing}
\quad \=\quad \=\quad \=\quad \=\quad \=\quad \=\quad\=\quad \kill
{$s \leftarrow 1$}\\
{$d \leftarrow true$}\\
\keyw{for} {$i \leftarrow n-1$} \keyw{to} {$0$} \keyw{do} \\
\> \keyw{if} {$(p_i = 1 ~xor~ d) = true$} \keyw{then}\\
\> \> {$d \leftarrow true$}\\
\> \keyw{else}\\
\> \> {$d \leftarrow false$}\\
\> \> {$s \leftarrow s +  2^i$}\\
end \keyw{for}\\
\keyw{return} {$s$}.
\end{tabbing}




\end{document}